\documentclass{article}
\usepackage{array,amsmath}
\numberwithin{equation}{section}
\usepackage{indentfirst}

\begin{document}

\newtheorem{thm}{Theorem}[section]
\newtheorem{cor}[thm]{Corollary}
\newtheorem{prop}[thm]{Proposition}
\newtheorem{conj}[thm]{Conjecture}
\newtheorem{lem}[thm]{Lemma}
\newtheorem{Def}[thm]{Definition}
\newtheorem{rem}[thm]{Remark}
\newtheorem{prob}[thm]{Problem}
\newtheorem{ex}{Example}[section]

\newcommand{\be}{\begin{equation}}
\newcommand{\ee}{\end{equation}}
\newcommand{\ben}{\begin{enumerate}}
\newcommand{\een}{\end{enumerate}}
\newcommand{\beq}{\begin{eqnarray}}
\newcommand{\eeq}{\end{eqnarray}}
\newcommand{\beqn}{\begin{eqnarray*}}
\newcommand{\eeqn}{\end{eqnarray*}}
\newcommand{\bei}{\begin{itemize}}
\newcommand{\eei}{\end{itemize}}

\newcommand{\pa}{{\partial}}
\newcommand{\V}{{\rm V}}
\newcommand{\R}{{\bf R}}
\newcommand{\K}{{\rm K}}
\newcommand{\e}{{\epsilon}}
\newcommand{\tomega}{\tilde{\omega}}
\newcommand{\tOmega}{\tilde{Omega}}
\newcommand{\tR}{\tilde{R}}
\newcommand{\tB}{\tilde{B}}
\newcommand{\tGamma}{\tilde{\Gamma}}
\newcommand{\fa}{f_{\alpha}}
\newcommand{\fb}{f_{\beta}}
\newcommand{\faa}{f_{\alpha\alpha}}
\newcommand{\faaa}{f_{\alpha\alpha\alpha}}
\newcommand{\fab}{f_{\alpha\beta}}
\newcommand{\fabb}{f_{\alpha\beta\beta}}
\newcommand{\fbb}{f_{\beta\beta}}
\newcommand{\fbbb}{f_{\beta\beta\beta}}
\newcommand{\faab}{f_{\alpha\alpha\beta}}

\newcommand{\pxi}{ {\pa \over \pa x^i}}
\newcommand{\pxj}{ {\pa \over \pa x^j}}
\newcommand{\pxk}{ {\pa \over \pa x^k}}
\newcommand{\pyi}{ {\pa \over \pa y^i}}
\newcommand{\pyj}{ {\pa \over \pa y^j}}
\newcommand{\pyk}{ {\pa \over \pa y^k}}
\newcommand{\dxi}{{\delta \over \delta x^i}}
\newcommand{\dxj}{{\delta \over \delta x^j}}
\newcommand{\dxk}{{\delta \over \delta x^k}}

\newcommand{\px}{{\pa \over \pa x}}
\newcommand{\py}{{\pa \over \pa y}}
\newcommand{\pt}{{\pa \over \pa t}}
\newcommand{\ps}{{\pa \over \pa s}}
\newcommand{\pvi}{{\pa \over \pa v^i}}
\newcommand{\ty}{\tilde{y}}
\newcommand{\bGamma}{\bar{\Gamma}}

\font\BBb=msbm10 at 12pt
\newcommand{\Bbb}[1]{\mbox{\BBb #1}}

\newcommand{\qed}{\hspace*{\fill}Q.E.D.}  

\title{Some inequalities on Finsler manifolds with weighted Ricci curvature bounded below \footnote{The first author is supported by the National Natural Science Foundation of China (11871126) and the Science Foundation of Chongqing Normal University (17XLB022) and the second author is supported by the National Natural Science Foundation of China (12071423)}}
\author{ Xinyue Cheng  \ \&  Zhongmin Shen}
\date{}

\maketitle

\begin{abstract}
We establish some important inequalities under a lower  weighted Ricci curvature bound on Finsler manifolds. Firstly, we establish a relative volume comparison of Bishop-Gromov type. As one of the applications, we obtain an upper bound for volumes of the Finsler manifolds. Further, when the S-curvature is bounded on the whole manifold, we obtain  a theorem of Bonnet-Myers type on Finsler manifolds. Finally, we obtain a sharp Poincar\'{e}-Lichnerowicz inequality by using integrated Bochner inequality, from which we obtain a sharp lower bound for the first eigenvalue on the Finsler manifolds. \\
{\bf Keywords:} Finsler metric; Ricci curvature; weighted Ricci curvature;  geodesic ball; volume comparison; Poincar\'{e}-Lichnerowicz inequality \\
{\bf MR(2010) Subject Classification:}  53B40, 53C60, 58C35
\end{abstract}

\section{Introduction}

Similar to Riemannian case, Finsler manifolds with Ricci curvature bounded below are always of some amazing properties. For example, the second author proved the Gromov-Bishop volume comparison theorem in Finsler geometry under the condition that ${\rm Ric} \geq (n-1)\lambda$. As an application, he obtained some precompactness and finiteness theorems for Finsler manifolds (\cite{shen}). Later, Ohta proved a new version of the Gromov-Bishop volume comparison theorem on Finsler manifolds with weighted Ricci curvature bounded below, ${\rm Ric}^{N} \geq K$ for some $K\in R$ and $N\in [ n, \infty)$ (\cite{Oh1}). On the other hand, by using Bochner-Weitzenb\"{o}ck type formula, Wang-Xia obtained a sharp lower bound for the first (nonzero) Neumann eigenvalue of Finsler-Laplacian in Finsler manifolds with weighted Ricci curvature ${\rm Ric}^{N}$ bounded below (\cite{WangXia}).  Q. Xia gave an (integrated) $p(>1)$-Bochner-Weitzenb\"{o}ck formula and the $p$-Reilly type formula on Finsler manifolds. As applications, she obtained the $p$-Poincar\'{e} inequality on an $n$-dimensional compact Finsler manifold without boundary or with convex boundary under the assumption that ${\rm Ric}^N \geq K$ for $N \in[n, \infty]$ and $K \in {R}$ (see \cite{Xia}).

In this paper, our main aim is to educe some important inequalities on Finsler manifolds with weighted Ricci curvature ${\rm Ric}^{\infty}$ bounded below. The paper actually consists of two parts: the first part is dedicated to the geometry of Finsler manifolds and the second part is dedicated to the analysis on Finsler manifolds.

Let $(M, F, m)$  be an $n$-dimensional  Finsler manifold with a smooth measure $m$ and $dm_{BH}$ denotes the Busemann-Hausdorff volume form of $F$. Then there is a positive smooth function $\phi$ on $M$ such that
\[
dm = \phi(x) dm_{BH}.
\]
Now, let
$ 0\leq \rho_o < t_o \leq +\infty$ and $\chi: (\rho_o, t_o)\to (0, +\infty)$ be a $C^{\infty}$ function  such that $\lim_{t\to  0^+} \chi(t)=0$ if $\rho_o =0$ and $\lim_{t\to t_o^{-}}\chi(t) =0 $ if $t_o < +\infty$.
Typical examples are $\chi(t)=\sin^{n-1}(t)$ or  $ \chi(t) =t^{n-1}$. Fix $p \in M$. Let $\rho (x):= d_{F}(p,x)$ be the distance function and $d_p:=\sup_{x\in M} d_{F}(p, x)$.  From now on, we always assume that $ d_p \leq t_o$ in our discussions. Then, under an upper bound on $\Delta \rho$ controlled by a function  $\chi(\rho)$  of $\rho$, we can obtain the following relative volume comparison of Bishop-Gromov type for the geodesic balls with respect to $m$.

\begin{thm} \label{revol} Let $(M, F, m)$ be an $n$-dimensional forward complete Finsler manifold  equipped with a smooth volume form $dm = \phi (x) dm_{BH}$ and $p\in M$.  Suppose that the distance function $\rho (x)=d_{F}(p, x)$ satisfies
\begin{equation}
\Delta \rho(x) \leq \frac{d}{d\rho} \Big [ \ln \chi(\rho)\Big ]|_{\rho =\rho(x)}, \ \  x\in \Omega_p\setminus B_p(\rho_o)  \label{Deltar}
\end{equation}
for a function $\chi(t)$ as above, where $\Omega_p$ denotes the cut-domain at $p$. Then for $\rho_o < r < R <  t_o$,
\be
\frac{{\rm Vol}(B_{p}(R)\setminus B_p(\rho_o))}{{\rm Vol} (B_{p}(r)\setminus B_p(\rho_o)) }  \leq \frac{\int_{\rho_o}^R \chi(t)dt}{\int_{\rho_o}^r \chi(t)dt}.  \label{Volume}
\ee
Moreover, if $\rho_o=0$ and  $\chi(t) = t^{n-1} \{ 1+ O(t)\}$, then  for $ 0 < r < R < t_o$
\be
{\rm Vol}(B_{p}(R)\setminus B_p(r)) \leq  \phi(p) \omega_{n-1} \int_r^R \chi(t)dt, \label{chit}
\ee
where $B_{p}(R)$ denotes the geodesic ball of radius $R$ at the center $p$ and $\omega_{n-1}={\rm Vol}(S^{n-1}(1))$ denotes the volume of the standard unit sphere in $R^n$.
\end{thm}

We must point out that, under certain curvature bounds on the Ricci curvature and the S-curvature, one can find $\chi(t) $ such that (\ref{Deltar}) holds (e.g. \cite{OS1}\cite{WeiW}\cite{WuXin}). However, in order to get an absolute upper bound on the volume of geodesic ball, one must get a function $\chi(t)$ with the property $\chi(t)=  t^{n-1} \{ 1+ O(t)\}$.

\begin{thm}\label{thm1.3}  Let $(M, F, m)$ be an $n$-dimensional forward complete  Finsler manifold  equipped with a smooth volume form $dm = \phi (x) dm_{BH}$ and $p\in M$.  Suppose that there are constants $ K >0$ and $\delta \geq 0$ such that
\be
 {\rm Ric}^{\infty} \geq K >0, \ \ \ \ {\bf S}|_{B_p(\frac{\pi}{2}\sqrt{\frac{n-1}{K}})} \geq -\delta ,\label{RV}
\ee
where ${\bf S}$ denotes the S-curvature of $F$. Then the volume of $M$ is  bounded by a constant depending on $n$, $K$ and $\delta$,
\[ {\rm Vol}(M) \leq \phi(p) K^{-n/2} c\Big (n, \frac{\delta}{\sqrt{K}} \Big ),\]
where $c\Big (n, \frac{\delta}{\sqrt{K}} \Big )$ is a constant depending only on $n$ and $\frac{\delta}{\sqrt{K}}$.
\end{thm}

It is known that  for a complete Riemannian manifold with a volume form, if ${\rm Ric}^{\infty}$ is strictly positive, then the volume is finite. Theorem \ref{thm1.3} shows that this is true for Finsler manifolds  and gives an estimate of upper bound of the volume. However, the Ricci curvature bound ${\rm Ric}^{\infty} \geq K >0$ does not control the size of the  Finsler manifold and
the volume of  metric balls. Thus the lower bound on the S-curvature  cannot  dropped. Further,  if one assumes that $|{\bf S} |\leq \delta $ on the whole manifold, then $M$ must be compact with ${\rm Diam}(M) \leq C(n, K, \delta)$ and the volume is bounded too,  ${\rm Vol}(M) \leq \phi (p) K^{- n /2}C(n, \frac{\delta}{\sqrt{K}} )$.  See Theorem \ref{prop4.5} below for details.

\vskip 2mm

Finally,  let us consider an analytic problem on Finsler manifolds  with weighted Ricci curvature ${\rm Ric}^{\infty}$ bounded below.  Starting from (\ref{Intform}), we can get the following inequality.

\begin{thm}\label{PLineq} {\rm (Poincar\'{e}-Lichnerowicz inequality)} Suppose that $(M, F, m)$ is closed and satisfies $m(M)=1$ and $\operatorname{Ric}^{\infty} \geq K>0$. Then, for any $f \in H^{1}(M)$, we have the following
\be
\operatorname{Var}_{m}(f) \leq \frac{1}{K} \int_{M} F^{2}({\nabla} f) d{m}-\frac{2}{K}\int_{M}\left(\int_{0}^{\infty} g(t)dt\right)dm, \label{PoLiin}
\ee
where ${\rm Var}_{m}(f)$ denotes the variance of $f$ and
\[
g(t):= g_{\nabla u_{t}}\left(\nabla ^{{\nabla u_{t}}}F({\nabla u_{t}}), \nabla ^{{\nabla u_{t}}}F({\nabla u_{t}})\right),
\]
$\left(u_{t}\right)_{t \geq 0}$ is the global solution to the heat equation with $u_{0}=f$.
\end{thm}

In \cite{Oh2}, Ohta obtained Poincar\'{e}-Lichnerowicz inequality under the condition that ${\rm Ric}^{\infty} \geq K>0$  as follows
\be
{\rm Var}_{m}(f) \leq \frac{1}{K} \int_{M} F^{2}({\nabla} f) d{m}. \label{eq15.8}
\ee
Note that (\ref{PoLiin}) is stronger than (\ref{eq15.8}).

For a closed Finsler manifold $(M, F, m)$ with $m (M) =1$, the first eigenvalue $\lambda_1$ is defined by (\cite{Cheng1}\cite{shen2})
\[ \lambda_1 := \inf_{f} \frac{\int_M F^{2}({\nabla} f) d{m}}{\operatorname{Var}_{m}(f)}.\]
Under the assumption that ${\rm Ric}^{\infty} \geq K >0$, we get from (\ref{eq15.8}) that
\[  \lambda_1 \geq K.\]
However, according to Theorem \ref{PLineq},  the above estimate is not sharp, namely, we can get from (\ref{PoLiin}) that $\lambda_1 \geq K+\delta $,  where
\[ \delta := \inf_{f} \frac{2}{{\rm Var}_m(f)}\int_{M}\left(\int_{0}^{\infty} g(t)dt\right)dm.\]

\vskip 3mm

The paper is organized as follows. In Section \ref{Sec2}, we give some necessary definitions and notations. Then we derive a relative volume comparison of Bishop-Gromov type and give the proof of Theorem \ref{thm1.3} based on Theorem \ref{revol} in Section \ref{Comparison}. Further, when the S-curvature is bounded on the whole manifold, we obtain  a theorem of Bonnet-Myers type on Finsler manifolds (that is, Theorem \ref{prop4.5}) in Section \ref{BMth}. Finally, the proof of a sharp Poincar\'{e}-Lichnerowicz inequality and a sharp lower bound for the first eigenvalue on Finsler manifolds are given in Section \ref{Sec3}.

\section{Preliminaries}\label{Sec2}

Let $M$ be an $n$-dimensional manifold. A Finsler metric $F$ on $M$ is a non-negative function on $TM$ such that $F$ is ${\cal C}^{\infty}$ on $TM\backslash \{0\}$ and the restriction $F_{x}:=F|_{T_{x}M}$ is a Minkowski function on $T_{x}M$ for all $x\in M$. For Finsler metric $F$ on $M$, there is a Finsler co-metric $F^{*}$ on $M$ which is non-negative function on the cotangent bundle $T^{*}M$ given by
\be
F^{*}(x, \xi):=\sup\limits_{y\in T_{x}M\setminus \{0\}} \frac{\xi (y)}{F(x,y)}, \ \ \forall \xi \in T^{*}_{x}M.
\ee
We call $F^{*}$ the dual Finsler metric of $F$. Finsler metric $F$ and its dual Finsler metric $F^{*}$ satisfy the following relation.

\begin{lem}{\rm (Lemma 3.1.1, \cite{shen2})}\label{shen311} Let $F$ be a Finsler metric on $M$ and $F^{*}$ its dual Finsler metric. For any vector $y\in T_{x}M\setminus \{0\}$, $x\in M$, the covector $\xi =g_{y}(y, \cdot)\in T^{*}_{x}M$ satisfies
\be
F(x,y)=F^{*}(x, \xi)=\frac{\xi (y)}{F(x,y)}. \label{shenF311}
\ee
Conversely, for any covector $\xi \in T_{x}^{*}M\setminus \{0\}$, there exists a unique vector $y\in T_{x}M\setminus \{0\}$ such that $\xi =g_{y}(y, \cdot)\in T^{*}_{x}M$.
\end{lem}

Naturally, by Lemma \ref{shen311}, we define a map ${\cal L}: TM \rightarrow T^{*}M$ by
\[
{\cal L}(y):=\left\{
\begin{array}{ll}
g_{y}(y, \cdot), & y\neq 0, \\
0, & y=0.
\end{array} \right.
\]
It follows from (\ref{shenF311}) that
\be
F(x,y)=F^{*}(x, {\cal L}(y)).
\ee
Thus ${\cal L}$ is a norm-preserving transformation. We call ${\cal L}$ the Legendre transformation on Finsler manifold $(M, F)$.

Take a basis $\{{\bf b}_{i}\}^{n}_{i=1}$ for $TM$ and its dual basis $\{\theta ^{i}\}_{i=1}^{n}$ for $T^{*}M$. Express $\xi ={\cal L}(y)=\xi _{i}\theta ^{i}$. Then
\be
\xi _{i} =g_{ij}(x,y)y^{j},
\ee
where $g_{ij}(x,y):=\frac{1}{2}\left[F^2\right]_{y^{i}y^{j}}(x,y)$. Let
\be
g^{*kl}(x,\xi):=\frac{1}{2}\left[F^{*2}\right]_{\xi _{k}\xi_{l}}(x,\xi).
\ee
For any $\xi ={\cal L}(y)$, differentiating $F^{2}(x,y)=F^{*2}(x,{\cal L}(y))$ with respect to $y^i$ yields
\[
\frac{1}{2}\left[F^2\right]_{y^{i}}(x,y)=\frac{1}{2}\left[F^{*2}\right]_{\xi _{k}}(x,\xi)g_{ik}(x,y), \label{shen3121}
\]
which implies
\be
g^{*kl}(x,\xi)\xi _{l}=\frac{1}{2}\left[F^{*2}\right]_{\xi _{k}}(x,\xi)=\frac{1}{2}g^{ik}(x,y)\left[F^2\right]_{y^{i}}(x,y)=y^{k}. \label{dualY}
\ee
Then,  we can get (see \cite{Cheng1}\cite{shen2})
\be
g^{*kl}(x,\xi)= g^{kl}(x,y), \label{Fdual}
\ee
where $(g^{kl}(x,y))= (g_{kl}(x,y))^{-1}.$

Given a smooth function $u$ on $M$, the differential $du_{x}$ at any point $x\in M$,
\[
du_{x}=\frac{\pa u}{\pa x^i}(x)dx^{i}
\]
is a linear function on $T_{x}M$. We define the gradient vector $\nabla u(x)$ of $u$ at $x\in M$ by $\nabla u(x):={\cal L}^{-1}\left(du(x)\right)\in T_{x}M$. In a local coordinate system, by (\ref{dualY}), we can express $\nabla u$ as
\be
\nabla u(x)=\left\{
\begin{array}{ll}
g^{*ij}(x,d u)\frac{\pa u}{\pa x^i}\frac{\pa}{\pa x^j}, & x\in M_{u},  \label{gradientV}\\
0, & x \in M\setminus M_{u},
\end{array}
\right.
\ee
where $M_{u}=\{x\in M \mid du (x)\neq 0\}.$ Further, by Lemma \ref{shen311}, we have the following
\be
du_{x}(v)=g_{\nabla u_{x}}(\nabla u_{x}, v), \ \ \ \forall v\in T_{x}M
\ee
and
\be
F(x, \nabla u_{x})=F^{*}(x, du_{x})=\frac{du_{x} (\nabla u_{x})}{F(x,\nabla u_{x})}.
\ee
We must be careful when $d u(x)=0,$ because $g^{*i j}(d u(x))$ is not defined and the Legendre transform $\mathcal{L}^{-1}$ is only continuous at the zero section. Besides, if $u \in \mathcal{C}^{l}(M),$ then ${\nabla}u$ is $\mathcal{C}^{l-1}$ on $M_{u}$ whereas only continuous on $M \backslash M_{u}$.

Let $(M, F, {m})$ be a Finsler manifold equipped with a measure ${m}$ on $M.$  Given an open set $\Omega \subset M,$ let $H_{\rm loc }^{1}(\Omega)$ be the space of weakly differentiable functions $u$ on $\Omega$ such that both $u$ and $F^{*}(d u)$ belong to $L_{\rm loc }^{2}(\Omega)$.  Define the energy functional $\mathcal{E}_{\Omega}: H_{\rm loc }^{1}(\Omega) \longrightarrow[0, \infty]$ by
\be
\mathcal{E}_{\Omega}(u):=\frac{\int_{\Omega} F^{*}(x, d u)^{2} d {m}}{\int_{\Omega} u^{2} dm}= \frac{\int_{\Omega} F(\nabla u)^{2} d {m}}{\int_{\Omega} u^{2} dm}.
\ee
Further, define the Sobolev space $H^{1}(\Omega):= W^{1, 2}(\Omega)$ as a space of functions $u \in  L^{2}(M)$ with $\int_{M}[F(\nabla u)]^{2} d m< \infty$, which is a Banach space with respect to the (absolutely homogeneous) Sobolev norm
$$
\|u\|_{H^{1}(\Omega)}:= \left\{\|u\|_{L^{2}(\Omega)}+ \frac{1}{2}\|F(\nabla u)\|_{L^{2}(\Omega)}+\frac{1}{2}\|F\left(\nabla (- u)\right)\|_{L^{2}(\Omega)} \right\}^{1 / 2}.
$$
Let $H_{0}^{1}(\Omega)$ be the closure of $\mathcal{C}_{c}^{\infty}(\Omega)$ in $H^{1}(\Omega)$,
where $\mathcal{C}_{c}^{\infty}(\Omega)$ denotes the set of $\mathcal{C}^{\infty}$-functions on $\Omega$ with compact support.

Associated with the measure $m$ on $M$, we first decomposed the volume form $d m$ of  $m$ as $d m = {\rm e}^{\Phi} dx^{1} dx^{2} \cdots d x^{n}$. Then the divergence of a differentiable vector field $V$ on $M$ is defined by
\be
{\rm div}_{m} V:= \frac{\partial V^{i}}{\partial x^{i}}+V^{i} \frac{\partial \Phi}{\partial x^{i}}, \quad V= V^{i} \frac{\partial}{\partial x^{i}}. \label{divdef1}
\ee
One can also define ${\rm div}_{m} V$ in the weak form by following divergence formula:
\be
\int_{M} \phi \ {\rm div}_{m}V \ d{m}=-\int_{M} d\phi (V) \ d{m} \label{divdef2}
\ee
for all $\phi \in \mathcal{C}_{c}^{\infty}(M)$.

Now we define the Finsler Laplacian $\Delta u$ of $u \in H_{\rm loc}^{1}(M)$ by
\be
\Delta u:= {\rm div}_{m}(\nabla u). \label{laplace1}
\ee
Equivalently, we can define Laplacian $\Delta u$ on the whole $M$ in the weak sense by
\be
\int_{M} \phi \ {\Delta}u \ d{m}:=-\int_{M} d \phi({\nabla} u) d{m} \label{laplace2}
\ee
for all $\phi \in \mathcal{C}_{c}^{\infty}(M)$.

From (\ref{laplace1}), Finsler Laplacian is a nonlinear elliptic differential operator of the second order. Moreover, since the gradient vector field $\nabla u$ is merely continuous on $M \backslash M_{u}$, even when $u \in \mathcal{C}^{\infty}(M),$ it is necessary to introduce the Laplacian in the weak form as (\ref{laplace2}).

By the definitions, we have the following for any smooth function $\varphi$ on $M$,
\be
{\rm div}_{m}(\varphi \nabla u)=\varphi \Delta u+d\varphi (\nabla u). \label{DLOper}
\ee

Let $V=V^{i} \frac{\partial}{\partial x^{k}}$ be a nonzero measurable vector field on $M.$ One can introduce a weighted Riemannian metric $g_{V}$ on $M$ via
\[
g_{V}(X, Y)=g_{i j}(x, V) X^{i} Y^{j}, \quad \text{ for } X, Y \in T_{x} M.
\]
In particular, $g_{V}(V, V)=F^{2}(x, V)$.  For $u \in H_{\rm loc }^{1}(M)$ such that $du=0$ almost everywhere on $\{x \in M \mid V(x)=0\}$ (in other words, $V \neq 0$ almost everywhere on $M_{u}$ ), we can define the linearized gradient vector and the linearized Laplacian on the weighted Riemannian manifold $\left(M, g_{V}, {m}\right)$ by
\be
\nabla^{V} u:=\left\{\begin{array}{ll}  g^{i j}(x, V) \frac{\partial u}{\partial x^{j}} \frac{\partial}{\partial x^{i}} & \text{on } M_{u}, \\
0 & \text{on } M \backslash M_{u},\end{array} \quad \Delta^{V} u:= {\rm div}_{m}\left(\nabla^{V} u\right)\right.. \label{weigraLa}
\ee

It is easy to get the following by (\ref{weigraLa}): for any $u \in H_{\rm loc}^{1}(M)$, we have
\be
\nabla^{\nabla u} u=\nabla u , \ \ \ \Delta^{\nabla u} u=\Delta u. \label{weiGra}
\ee
Further, for any $f_{1}, f_{2} \in H_{\mathrm{loc}}^{1}(M)$ satisfying $d f_{1}=d f_{2}=0$ almost everywhere on $M \backslash M_{u}$, we have
\be
d f_{2}\left(\nabla^{\nabla u} f_{1}\right)= d f_{1}\left(\nabla^{\nabla u} f_{2}\right). \label{LinGra}
\ee

\vskip 3mm

Finsler geometry is just Riemannian geometry without the quadratic restriction (\cite{Chern}). The Ricci curvature in Finsler geometry is just a natural extension of the Ricci curvature in Riemann geometry. However, we have a difficulty on the choice of a measure in Finsler geometry, because it is impossible to choose a unique canonical measure like the volume measure in the Riemannian setting. Naturally, by choosing an arbitrary measure $m$ on a Finsler manifold, Ohta modified the Ricci curvature and defined the weighted Ricci curvature in Finsler geometry (\cite{Oh1}). Concretely, for an $n$-dimensional Finsler manifold $(M, F, m)$ equipped with a smooth measure $m$ and for any $v \in T_{x} M \backslash\{0\},$ let $\eta:(-\varepsilon, \varepsilon) \longrightarrow M$ be the geodesic with $\dot{\eta}(0)=v$ and decompose the volume form $d m$ of measure ${m}$ along $\eta$ as
\[
dm = {\rm e}^{-\psi_{\eta}} \sqrt{{\rm det}\left(g_{i j}(\eta , \dot{\eta})\right)}\ d x^{1} d x^{2} \cdots d x^{n},
\]
where $\psi_{\eta}= \psi_{\eta}\left(\eta (t), \dot{\eta}(t)\right): (-\varepsilon, \varepsilon) \longrightarrow R$ is a $\mathcal{C}^{\infty}$-function. Then, for $N \in R \backslash\{n\},$ define the weighted Ricci curvature
\be
{\rm Ric}^N(v):= {\rm Ric}(v)+\psi_{\eta}^{\prime \prime}(0)-\frac{\psi_{\eta}^{\prime}(0)^{2}}{N-n}. \label{weiRicci1}
\ee
As the limits of $N \rightarrow \infty$ and $N \downarrow n,$  we define the weighted Ricci curvatures as follows.
\be
\begin{aligned}
{\rm Ric}^{\infty}(v) &: = {\rm Ric}(v)+\psi_{\eta}^{\prime \prime}(0), \\
{\rm Ric}^n(v)& :=\left\{\begin{array}{ll} {\rm Ric}(v)+\psi_{\eta}^{\prime \prime}(0) & \text { if } \psi_{\eta}^{\prime}(0)=0, \\
-\infty & \text { if } \psi_{\eta}^{\prime}(0) \neq 0. \end{array}\right.
\end{aligned}
\ee

Note that the quantity $\psi_{\eta}^{\prime}(0)={\bf S}(x, v)$ is just the S-curvature with respect to the measure $m$ and $\psi_{\eta}^{\prime \prime}(0)= \dot{\bf S}(x, v)= {\bf S}_{| m}(x, v)v^{m}$ in Finsler geometry, where $`` | "$ denotes the horizontal covariant derivative with respect to the Chern connection  (\cite{ChernShen}\cite{Oh1}). Hence we can rewrite ${\rm Ric}^N(v)$ as
\be
{\rm Ric}^N(v):= {\rm Ric}(v)+\dot{\bf S}(x, v) -\frac{{\bf S}(x, v)^{2}}{N-n}. \label{weiRicci2}
\ee

From (\ref{weiRicci2}), we can get the following further description about the weighted Ricci curvature in Finsler geometry. Let $(M, F, m)$ be an $n$-dimensional Finsler manifold with  $d m= \sigma (x) dx^{1} \cdots dx^{n}$.
Let $Y$ be a $C^{\infty}$ geodesic field on an open subset $U \subset M$ and $\hat{g}=g_{Y} .$  Let
$$
d m:=e^{- \psi} {\rm Vol}_{\hat{g}}, \ \ \ {\rm Vol}_{\hat{g}}= \sqrt{{det}\left(g_{i j}\left(x, Y_{x}\right)\right)}dx^{1} \cdots dx^{n}.
$$
It is easy to see that $\psi$ is given by
$$
\psi (x)= \ln \frac{\sqrt{\operatorname{det}\left(g_{i j}\left(x, Y_{x}\right)\right)}}{\sigma(x)}=\tau\left(x, Y_{x}\right),
$$
which is just the distortion along $Y_{x}$ at $x\in M$ (\cite{ChernShen}).
Let $y := Y_{x}\in T_{x}M$ (that is, $Y$ is a geodesic extension of $y\in T_{x}M$). Then, by the definitions of the S-curvature and the Hessian (\cite{shen1}\cite{shen2}), we have
\beqn
&& {\bf S}(x, y)=y[\tau(x, Y_{x})]= d \psi (y),\\
&& \dot{\bf S}(x, y)=y[{\bf S}(x, Y)]=y[Y(\psi)]= {\rm Hess} \psi (y).
\eeqn
Hence,
\be
{\rm Ric}^N(y)= {\rm Ric}(y)+ {\rm Hess}  \psi(y)-\frac{d \psi(y)^{2}}{N-n}.   \label{weRicci3}
\ee
Obviously, (\ref{weRicci3}) is an analogue of the weighted Ricci curvature in Riemannian geometry.

\section{A relative volume comparison of Bishop-Gromov type and its applications}\label{Comparison}

In this section,  firstly we will derive a relative volume comparison of Bishop-Gromov type and prove Theorem \ref{revol}.

Let $(M, F, m)$ be an $n$-dimensional forward complete Finsler manifold with a smooth measure $m$.  Let $dm_{BH}$ denote the volume form of the Busemann-Hausdorff measure  $m_{BH}$ of $F$. Then the volume form $dm$ of $m$ is a multiple of $dm_{BH}$,
\[  dm = \phi (x) dm_{BH}.\]
These measures induce measures on hypersurfaces if a normal vector field is chosen. Let $p\in M$. The volumes of a small geodesic sphere $S_p(r)$ and a small geodesic ball $B_p(r)$ with respect to $m$ are given by
\begin{eqnarray}
{\rm Vol}(S_p(r)) & = & \phi(p) \omega_{n-1} r^{n-1}\Big \{ 1+ O(r)\Big \}, \label{v1}\\
 {\rm Vol} (B_p(r)) & =  &\phi(p)\frac{\omega_{n-1}}{n} r^n \Big \{ 1 +O(r)\Big \}, \label{v2}
\end{eqnarray}
where $\omega_{n-1}={\rm Vol}(S^{n-1}(1))$ and $S_p(r)$ denotes the geodesic  sphere of radius $r$ at center $p$.  See  \cite{shen2}.

Let  $\Omega_p$ be the cut-domain -- the open domain enclosed by the cut-locus ${\rm Cut}(p)$ at $p$.  The distance function $\rho (x)=d_F(p, x)$ is $C^{\infty}$ on $\Omega_p\setminus\{p\}$.
Then $\tilde{S}_p(t) := S_p( t)\cap \Omega_p$ is a $C^{\infty}$ hypersurface in $M$ for $t >0$. The  volume form $dm$ induces a volume form $dA_t$ on $\tilde{S}_p(t)$. On the other hand, the volume form $dm$ determines the  volume form $dm_p$ on $T_{p}M$. $dm_p$ induces  the volume form $dA_{p}$ on $S_pM := \{ y \in T_{p}M ~| \ F(p, y)=1\}$.
Define a map $\varphi_t: S_pM \to M$  by
\be
 \varphi_t(y) := \exp_p(ty), \ \ \ \ \ y\in S_pM.
\ee
Set
\be
(\varphi_t)^* dA_t =  \eta_t(y) dA_p.
\ee
The function $\eta_t(y)$ is positive for $ 0 < t < i_y$ (the cut-value of $y$). For a small $ s >0$,
\[ {\rm Vol}(S_p( s)) =\int_{S_p( s)} dA_s = \int_{S_pM} \eta_s (y) dA_p.\]
Further, the Laplace of $\rho (x)$ with  respect to $F$ and $dm$ can be expressed by
\be
\Delta \rho |_{\tilde{S}_{p}(t)} =  \frac{d}{dt} \Big [ \ln \eta_t(y) \Big ].  \label{Deltaeta}
\ee
This is the link between the volume of geodesic spheres  $S_{p}(r)$  and the Laplace of the distance function  $\Delta \rho (x)$ (see (16.6) and Proposition 14.3.1 in \cite{shen2}).

Let  $0\leq \rho_o < t_o \leq \infty$ and  $\chi: (\rho_o, t_o)\to (0,\infty)$ be a $C^{\infty}$ function  such that   $\lim_{t\to 0^{+}} \chi(t)=0$ if $\rho_o=0$ and  $\lim_{t\to t_o^+} \chi(t)=0$ if $t_o < \infty$.
Assume that $ d_p:=\sup_{x\in M} d_{F}(p,  x) \leq t_o$ and
\be
\Delta \rho (x) \leq \frac{d}{dt} \Big [\ln \chi(t)\Big ]|_{t = \rho (x)}, \ \ \ \ \\ \forall x\in \Omega_p \setminus B_p(\rho_o).  \label{Deltachi}
\ee
By (\ref{Deltaeta}) and (\ref{Deltachi}), we get
\[ \frac{d}{dt} \Big [ \ln \eta_{t}(y)\Big ] \leq \frac{\chi'(t)}{\chi(t)}, \ \ \ \ \ \rho_o < t < i_y.\]
This implies that
\[  \frac{d}{dt} \Big [ \ln \frac{\eta_t(y)}{\chi(t)} \Big ] \leq 0, \ \ \ \ \ \rho_o < t < i_y.\]
Thus $\eta_t(y) /\chi(t)$ is non-increasing monotonically.  Let $\tilde{\eta}_t(y) :=\eta_t(y)$ for $\rho_o \leq t < i_y$ and $\tilde{\eta}_t(y) =0$ for $t \geq i_y$. Then  $\tilde{\eta}_t(y)/\chi(t)$ is non-increasing monotonically for $ \rho_o < t < t_o$, that is,
\[
\frac{\tilde{\eta}_t(y)}{\chi(t)} \leq \frac{\tilde{\eta}_s(y)}{\chi(s)}, \ \ \ \ \ \rho_o < s < t < t_o.
\]
Integrating it over $S_pM$ with respect to $dA_p$, we get
\be
\frac{ {\rm Vol}(\tilde{S}_p(t))}{\chi(t)} \leq \frac{{\rm Vol}(\tilde{S}_p(s))} {\chi(s)}, \ \ \ \ \ \rho_o < s < t < t_o. \label{Volchi}
\ee
Rewrite (\ref{Volchi}) as
\be
{\rm Vol}(\tilde{S}_p(t)) \chi(s) \leq  {\rm Vol}(\tilde{S}_p(s)) \chi(t), \ \ \ \ \ \ \rho_o < s < t < t_o. \label{Volchi1}
\ee
Let $ \rho_o < r < R < t_o$.  Fix $ t\in [r , R]$. Integrating (\ref{Volchi1}) with respect to $ s \in [\rho_o, r]$, we get
\[
{\rm Vol}(\tilde{S}_p(t)) \int_{\rho_o}^r \chi(s)ds \leq {\rm Vol}(B_{p}(r)\setminus B_p(\rho_o)) \chi(t), \ \ \ \ \ \ r \leq t \leq R.
\]
Then integrating it with respect to $ t \in [r, R]$, we get
\be
{\rm Vol} (B_{p}(R)\setminus B_{p}(r)) \int_{\rho_o}^r \chi(s) ds \leq {\rm Vol}(B_{p}(r)\setminus B_p(\rho_o)) \int_r^R \chi(t) dt. \label{vv}
\ee
By adding ${\rm Vol}(B_{p}(r)\setminus B_p(\rho_o)) \int_{\rho_o}^r \chi(t) dt$ to both sides of (\ref{vv}), we get
\[
{\rm Vol} (B_{p}(R)\setminus B_p(\rho_o)) \int_{\rho_o}^r \chi(t) dt \leq {\rm Vol}(B_{p}(r)\setminus B_p(\rho_o)) \int_{\rho_o}^R \chi(t) dt.
\]
Then, we obtain
\[
 \frac{ {\rm Vol}(B_{p}(R)\setminus B_p(\rho_o))}{{\rm Vol}(B_{p}(r)\setminus B_p(\rho_o)) } \leq \frac{\int_{\rho _{0}}^{R}\chi(t)dt
 }{\int_{\rho _{0}}^{r}\chi(t)dt }.\]
It is just (\ref{Volume}).

From (\ref{vv}), we also get
\be
\frac{{\rm Vol} (B_{p}(R)\setminus B_{p}(r)) }{\int_r^R \chi(t) dt} \leq \frac{{\rm Vol}(B_{p}(r)\setminus B_p(\rho_o))}{\int_{\rho_o}^r \chi(t) dt}.  \label{vv0}
\ee

Now, assume that $\rho_o=0$. It follows from (\ref{vv0}) that
\be
\frac{{\rm Vol} (B_{p}(R)\setminus B_{p}(r)) }{\int_r^R \chi(t) dt} \leq \frac{{\rm Vol}(B_{p}(r))}{\int_{0}^r \chi(t) dt}.  \label{vv1}
\ee
Further, assume that
\[ \chi(t) = t^{n-1} \{ 1+ O(t)\Big \}.
\]
By (\ref{v2}), we have
\[
\lim_{r\to 0} \frac{ {\rm Vol}(B_{p}(r))}{\int_0^r\chi(t)dt } = \phi(p)\omega_{n-1}.
\]
Letting $r\to0^+$ in (\ref{vv1}) yields
\be
\frac{{\rm Vol} (B_{p}(R)) }{\int_0^R \chi(t) dt} \leq \phi(p) \omega_{n-1}.  \label{vv2}
\ee
Replacing $R$ by $r$ in (\ref{vv2}) and applying the result to (\ref{vv1}), we get
\be
\frac{{\rm Vol} (B_{p}(R)\setminus B_{p}(r)) }{\int_r^R \chi(t) dt} \leq \phi(p)\omega_{n-1}.  \label{vv3}
\ee
This gives rise to  (\ref{chit}).  This completes the proof of Theorem \ref{revol}.
\vskip 2mm

In the following, we will show that we can find some $\chi (t)$ satisfying the conditions in Theorem \ref{revol}. Firstly, define
\be
{s}_{c}(t):=\left\{\begin{array}{ll}
\frac{1}{\sqrt{c}} \sin (\sqrt{c} t) & \text { for } c >0, \\
t & \text { for } c =0, \\
\frac{1}{\sqrt{- c}} \sinh (\sqrt{- c} t) & \text { for } c <0 .
\end{array}\right.
\ee
Note that ${s}_{c}(t)$ is the solution to the differential equation
$$
f^{\prime \prime}+c f=0, \quad f(0)=0, \quad f^{\prime}(0)=1.
$$
Further, define
\be
ct_{c}(t):= \frac{{s}'_{c}(t)}{{s}_{c}(t)}= \left\{\begin{array}{ll}
\sqrt{c} \cdot \operatorname{cotan}(\sqrt{c} t)  & \text { for } c>0 , \\
\frac{1}{t}  & \text { for } c=0,  \\
\sqrt{-c} \cdot \operatorname{cotanh}(\sqrt{-c} t)  & \text { for } c<0 .
\end{array}\right.
\ee

Let $\hat{g}:=g_{\nabla \rho}$ denote the induced Riemannian metric on $\Omega_p$ and $\hat{\Delta}$ denote the Laplacian on $(M, \hat{g}, {\rm Vol}_{\hat{g}})$. Let $dm = e^{-f}{\rm Vol}_{\hat{g}}$.
It is proved in \cite{shen2} that  $\nabla\rho =\hat{\nabla}\rho$ is a geodesic field of $F$ and $\hat{g}$ and
\be
{\bf S}(x, \nabla \rho_x) = df_x(\hat{\nabla} \rho_x), \ \  {\rm Ric}(x, \nabla \rho_x) = \widehat{\rm Ric}(x, \hat{\nabla}\rho_x), \ \  {\rm Hess} f(\nabla \rho_x)= \widehat{\rm Hess} f (\hat{\nabla}\rho_x).
\ee
Then
\begin{eqnarray*}
{\rm Ric}^{\infty} (x, \nabla\rho_x) &  =  & \widehat{\rm Ric}^{\infty} (x, \hat{\nabla}\rho_x):=\widehat{\rm Ric}(x, \hat{\nabla}\rho_x)+ {\rm Hess} f(\nabla \rho_x)\\
{\rm Ric}^{N} (x, \nabla\rho_x) &  =  & \widehat{\rm Ric}^{N} (x, \hat{\nabla}\rho_x):=\widehat{\rm Ric}(x, \hat{\nabla}\rho_x)+ {\rm Hess} f(\nabla \rho_x) -\frac{1}{N-n} [ df_x(\hat{\nabla} \rho_x)]^2.
\end{eqnarray*}
Further,
\be
\Delta \rho =\hat{\Delta}_f \rho :=\hat{\Delta}\rho -df(\hat{\nabla}\rho)= \hat{\Delta}\rho - df_x(\hat{\nabla} \rho_x).
\ee

\bigskip

Under the lower  Ricci curvature  bound ${\rm Ric}\geq (n-1)c$,  we have $\widehat{\rm Ric}(x, \hat{\nabla}\rho_x) \geq (n-)c$.
Then
\[ \hat{\Delta} \rho \leq (n-1)\operatorname{ct}_{c}(\rho).\]
Under the lower S-curvature bound ${\bf S}\geq -\delta$, we have $df_x(\hat{\nabla} \rho_x)\geq -\delta$.  Then
\[ \Delta \rho \leq (n-1)\operatorname{ct}_{c}(\rho)+\delta.\]
This leads to the following

\begin{lem}\label{lem4.1} {\rm (Laplacian comparison, \cite{shen2}\cite{WuXin})}\label{WuXin}  \  Let $(M, F, dm)$ be an $n$-dimensional Finsler manifold with Ricci curvature satisfying ${\rm Ric} \geq(n-1) c $ and ${\bf S}\geq -\delta$.  Then the following holds whenever the distance function $\rho = \rho (x)$ is smooth
\be
\Delta \rho \leq \frac{d}{d t} \Big [ \ln \chi(t)\Big ]|_{t=\rho(x)}, \label{WuXin2}
\ee
where $\chi(t)=  s_c(t)^{n-1} e^{\delta t} $, $0 < t < d_p$.
\end{lem}

Then, by Theorem \ref{revol} and Lemma \ref{WuXin}, we can obtain a volume comparison theorem of Bishop-Gromov type given in \cite{shen}.

\bigskip
Assume  the lower Ricci curvature bound ${\rm Ric}^{\infty} \geq (n-1)c$,  then $\widehat{\rm Ric}^{\infty} (x, \hat{\nabla}\rho_x) \geq (n-1)c$. Then we obtain the following lemma by applying a result in \cite{WeiW}.

\begin{lem} \label{lem4.3} {\rm (Laplacian comparison, \cite{WeiW}\cite{Yin})} Let $(M, F, m)$ be an $n$-dimensional Finsler manifold with weighted Ricci curvature satisfying ${\rm Ric}^{\infty} \geq(n-1) c.$  Then the following bound on $\Delta \rho$ holds.
\ben
\item[(a)] If ${\bf S}\geq -\delta$, then the following holds on $\Omega_p \cap B_p(r_o)$,
\[ \Delta\rho \leq \frac{d}{dt} \Big [ \ln \chi(t)\Big ]|_{t=\rho(x)},\]
where $\chi(t) = [s_c(t)]^{n-1} e^{\delta t}$, $ 0 < t < r_o$. Here $ r_o=+\infty$ when $c \leq 0$ and $r_o = \frac{\pi}{2 \sqrt{c}}$ when $c >0$.
\item[(b)] If the distortion $\tau$ satisfies that $|\tau |\leq k$,  then the following holds on $\Omega_p\cap B_p(r_o)$,
\[
\Delta\rho \leq  \frac{d}{dt} \Big [ \ln \chi(t)\Big ]|_{t=\rho(x)},\]
where $\chi(t) := [s_c(t)]^{n+4k-1}$, $ 0 < t < r_o$.  Here $r_o:=+\infty$ when $ c \leq 0$ and $r_o=\pi/(4\sqrt{c})$ when $c > 0$.
\een
\end{lem}

\begin{rem}
Lemma \ref{lem4.3}(a) is similar to  Lemma \ref{lem4.1}.  However, in the case when $c >0$, it might be possible that $ r_o < d_p$. Thus Theorem \ref{revol} holds only for $ 0 \leq \rho_o < r < R < r_o$.  By Theorem  \ref{revol} and Lemma \ref{lem4.3}, we can obtain a volume comparison theorem of Bishop-Gromov type for metric balls of radius less than $r_o$.   See \cite{Yin}.
\end{rem}

\begin{lem}\label{lem4.4}{\rm (Laplacian comparison, \cite{OS1})}\label{OS1}  \ Let $(M, F, {m})$ be forward or backward complete and assume that ${\rm Ric}^N \geq (N-1) K$ for some $K \in {R}$ and $N \in[n, \infty).$ Then, for any $p \in M,$ the distance function $\rho (x):=d_{F}(p, x)$ satisfies the following on $\Omega_p$
\be
\Delta \rho (x) \leq \frac{d}{dt} \Big [ \ln \chi(t)\Big ]|_{t=\rho(x)}\label{Ohta}
\ee
where $\chi(t):= [s_K(t)]^{N-1}$, $0 < t < r_o$. Here $r_o=\pi/\sqrt{K}$ when $K >0$.
\end{lem}

 By Theorem \ref{revol} and Lemma \ref{lem4.4}, we can obtain a volume comparison theorem of Bishop-Gromov type given in \cite{Oh1}.

\bigskip

Assume that
\[ {\rm Ric}^{\infty} \geq K. \]
This implies that $\widehat{\rm Ric}^{\infty} (x, \hat{\nabla}\rho_x) \geq K$.
By Theorem 3.1 in \cite{WeiW},  along a minimal geodesic $\gamma(t)$, $ 0 \leq t \leq t_o$ from  $\gamma(0)=p$,  for any $0< \rho_o < r < t_o $,
\be
\hat{\Delta}_f\rho |_{\gamma(r)} \leq \hat{\Delta}_{f} \rho|_{\gamma(\rho_o)} - K (r-\rho_o).
\ee
Then
\[ \Delta \rho (x) \leq \Delta \rho (x')  - K\{ \rho(x)-\rho_o\}, \ \ \ \ \rho(x) > \rho_o,\]
where $x'\in S_p(\rho_o)$ is on the minimal geodesic from $p$ to $x$.
We assume that $ 0 < \rho_o < i_p$. Then  $S_p(\rho_o)$ is a smooth hypersurface in $\Omega_p$.
Let
\[ m_o := \sup_{\rho(x')=\rho_o} \Delta \rho (x') <+\infty.\]
Then
\[  \Delta \rho (x) \leq m_o -K\{ \rho(x) -\rho_o\}, \ \  \ \ \rho(x) > \rho_o.\]
This leads to the following

\begin{lem}\label{lem4.5} {\rm (Laplacian comparison, \cite{WeiW}\cite{Yin})}  \  Let $(M, F, m)$ be an $n$-dimensional Finsler manifold with weighted Ricci curvature satisfying ${\rm Ric}^{\infty} \geq K $.  Then, for any $p \in M$, the following holds whenever the distance function $\rho (x):= d_{F}(p, x)$ is smooth and $\rho(x) > \rho_o$
\be
\Delta \rho \leq \frac{d}{dt} \Big [ \ln \chi(t)\Big ]|_{t=\rho(x)}, \label{WuXin2}
\ee
where $\chi(t)=   e^{m_o (t-\rho_o) -\frac{1}{2}K (t-\rho_o)^2} $, $\rho_o < t < \infty$.
\end{lem}

\bigskip
\noindent
{\it Proof of Theorem \ref{thm1.3}}:
We assume that
\[ {\rm Ric}^{\infty} \geq K > 0, \ \ \ \ {\bf S}|_{B_p(\frac{\pi}{2}\sqrt{\frac{n-1}{K}} ) } \geq -\delta.\]
Then ${\rm Ric}^{\infty}\geq (n-1)c$  with $c = K/(n-1)$.

Let
\[  r:=\frac{\pi}{2}\sqrt{\frac{n-1}{K}}.\]
By Lemma \ref{lem4.3} (a), the following holds on $B_p(r )$,
\be
 \Delta \rho \leq  \frac{d}{dt} \Big [\ln \chi_o (t)\Big ] |_{t=\rho(x)}, \ \ \ \ \rho(x) < r\label{m_o}
\ee
where $\chi_o (t):= [s_c(t)]^{n-1} e^{\delta t}$ where $ c=\frac{K}{n-1}$.
Note that $\chi_o(t) = t^{n-1}\{ 1+O(t)\}$. By (\ref{vv2}),
\begin{eqnarray}
 {\rm Vol}(B_p(r)) &\leq  &\phi(p)\omega_{n-1} \int_0^{r} [s_c(t)]^{n-1}e^{\delta t}  dt\nonumber\\
& =  & \phi(p)\omega_{n-1}
\Big (  \frac{n-1}{K} \Big)^{n/2} \int_{0}^{\pi/2} \sin^{n-1}(s) e^{\delta \sqrt{\frac{n-1}{K}} s} ds,\label{small3}
\end{eqnarray}

Suppose that
\[  d_p \leq r.\]
Then
\[ {\rm Vol}(M) ={\rm Vol}(B_p(r)) \leq \phi(p) K^{-n/2} c\Big ( n, \frac{\delta}{\sqrt{K}}\Big ).\]

We now assume that
\[  d_p  > r=\frac{\pi}{2}\sqrt{\frac{n-1}{K}}.\]
In this case, the manifold might noncompact, i.e. $d_p=+\infty$.
Let $ \rho_o:= r/2$.
By (\ref{m_o}), we have
\[ m_o =\sup_{\rho(x) = \rho_o} \Delta \rho \leq\frac{d}{dt} \Big[ \ln \chi_o(t)\Big ]_{t=\rho_o}  =\sqrt{(n-1)K}+\delta.
\]
Let
\be
 \chi(t) := e^{(\sqrt{(n-1)K}+\delta)(t-\rho_o) -\frac{K}{2} (t-\rho_o)^2}.\label{chit1}
\ee
Then by Lemma \ref{lem4.5},  the following holds on $M\setminus B_p(\rho_o)$,
\be
\Delta \rho \leq \frac{d}{dt} \Big [ \ln \chi(t)\Big ]|_{t=\rho(x)},  \ \ \ \ \rho(x) > \rho_o.\label{WuXin2*}
\ee
It follows from
 Theorem \ref{revol} that for $ \rho _{0} < r < R \leq d_p$,
\be
{\rm Vol}(B_{p}(R)) \leq  {\rm Vol}(B_p(\rho_o))+\frac{{\rm Vol} (B_{p}(r)\setminus B_p(\rho_o))}{\int_{\rho_o}^r \chi(t)dt} {\int_{\rho_o}^R \chi(t)dt  }.  \label{Volume**}
\ee
It follows from (\ref{Volume**}) that
\be
{\rm Vol}(M) \leq  {\rm Vol}(B_p(\rho_o))+\frac{{\rm Vol} (B_{p}(r)\setminus B_p(\rho_o))}{\int_{\rho_o}^r \chi(t)dt} {\int_{\rho_o}^{\infty} \chi(t)dt }. \label{VVVV}
\ee

We now estimate each term in (\ref{VVVV}). Similar to (\ref{small3}), we have
\begin{eqnarray}
 {\rm Vol}(B_p(\rho_o)) &\leq  &\phi(p)\omega_{n-1} \int_0^{\rho_o} [s_c(t)]^{n-1}e^{\delta t}  dt\nonumber\\
& =  & \phi(p)\omega_{n-1} \Big (  \frac{n-1}{K} \Big)^{n /2} \int_{0}^{\pi/4} \sin^{n-1}(s) e^{ \delta\frac{\sqrt{n-1}}{\sqrt{K}} s} ds. \label{small1}
\end{eqnarray}
Further, by (\ref{chit}), we have
\begin{eqnarray}
 {\rm Vol}(B_p (r)\setminus B_p(\rho_o)) & \leq  &  \phi(p)\omega_{n-1} \int_{\rho _{0}}^{r} [s_c(t)]^{n-1}e^{\delta t} dt\nonumber \\
& = & \phi(p)\omega_{n-1}
\Big (  \frac{n-1}{K} \Big)^{n / 2} \int_{\pi/4}^{\pi/2} \sin^{n-1}(s) e^{\delta\frac{\sqrt{n-1}}{\sqrt{K}} s} ds. \label{small2}
\end{eqnarray}
By substitution, we have
\[ \int_{\rho_o}^r \chi(t) dt = \frac{1}{ \sqrt{(n-1)K} +\delta } \int_0^h e^{s-\frac{1}{2}K_os^2} ds,\]
\[ \int_{\rho_o}^{\infty} \chi(t) t = \frac{1}{ \sqrt{(n-1)K} +\delta } \int_0^{\infty} e^{s-\frac{1}{2}K_os^2} ds,\]
where $ h =\frac{\pi}{4} \sqrt{n-1}\left(\sqrt{(n-1)}+\delta/\sqrt{K} \right) $ and $K_o:= 1/( \sqrt{n-1}+\delta/\sqrt{K})^2$.
Then by (\ref{VVVV}), we get
\begin{eqnarray*}
{\rm Vol} (M) & \leq &  \phi(p)\omega_{n-1}
\Big (  \frac{n-1}{K} \Big)^{n/2}\Big \{ \int_{0}^{\pi/4} \sin^{n-1}(s) e^{\delta {\frac{\sqrt{n-1}}{\sqrt{K}}} s} ds \\
&& + \frac{\int_{\pi/4}^{\pi/2} \sin^{n-1}(s) e^{\delta {\frac{\sqrt{n-1}}{\sqrt{K}}} s} ds}{\int_0^h e^{s-\frac{1}{2}K_os^2} ds}    \int_0^{\infty} e^{s-\frac{1}{2}K_os^2} ds \Big \} \\
&=: & \phi(p)K^{-n/2} c\Big ( n, \frac{\delta}{\sqrt{K}} \Big ).
\end{eqnarray*}
This completes the proof of  Theorem \ref{thm1.3}.

\section{A Theorem of Bonnet-Myers type}\label{BMth}

Let $(M, F, m)$ be an $n$-dimensional forward complete Finsler manifold with a smooth volume form $dm = \phi (x) dm_{BH}$. The classical Bonnet-Myers Theorem in Finsler geometry says that, if ${\rm Ric}\geq (n-1)K >0, $ then $M$ is in fact compact and ${\rm Diam} (M) \leq \frac{\pi}{\sqrt{K}}$ (\cite{BCS}, Theorem 7.7.1).  If the weighted Ricci curvature satisfies ${\rm Ric} _{N} \geq (N-1)K>0, N \in [n, \infty),$ then ${\rm Diam}(M) \leq \frac{\pi}{\sqrt{K}}$ (\cite{Oh1}).

In this section, we will consider the case when the S-curvature is bounded on the whole manifold and give a theorem of Bonnet-Myers type.

\begin{thm}\label{prop4.5}
Let $(M, F, m)$ be an $n$-dimensional forward complete Finsler manifold with a smooth volume form $dm = \phi (x) dm_{BH}$ and $p\in M$. Suppose that
\[ {\rm Ric}^{\infty} \geq K >0, \ \ \ \ \ |{\bf S}|\leq \delta .\]
Then the diameter and the volume is bounded.
\be
{\rm Diam} (M) \leq \frac{\pi}{\sqrt{K}} \Big ( \frac{\delta}{\sqrt{K}} + \sqrt{ \frac{\delta^2}{K} + n-1}\Big), \label{diamM}
\ee
\be
{\rm Vol} (M) \leq  \phi(p)K^{-n/2} C\Big ( n, \frac{\delta}{\sqrt{K}} \Big ), \label{VolofM}
\ee
where $C\Big ( n, \frac{\delta}{\sqrt{K}} \Big )$ is a constant depending only on $n$ and $\frac{\delta}{\sqrt{K}}$.
\end{thm}
{\it Proof}:  Firstly let $N$ be an arbitrary number with
\[   N > n + \frac{\delta^2}{K}.\]
Let
\[   H:= \frac{ K(N-n) -\delta^2} { (N-1)(N-n)} >0.\]
Then
\be
  {\rm Ric}^N = {\rm Ric}^{\infty} - \frac{1}{N-n} {\bf S}^2
\geq K -\frac{1}{N-n} \delta^2 = (N-1)  H >0.\label{RicH}
\ee
By a theorem of Bonnet-Myers type proved by S. Ohta in \cite{Oh1},
\[ {\rm Diam} (M)\leq \frac{\pi}{\sqrt{H}}.\]
Then
\[ {\rm Diam}(M) \leq  f_{K, \delta} (N):= \pi \sqrt{\frac{(N-1)(N-n)}{K(N-n)-\delta^2}   }   .\]
The minimum value of $f_{K,\delta} (N) $ over $ n +\delta^2/K < N < +\infty$ is
\[  \min f_{K,\delta} = \frac{\pi\gamma}{\sqrt{K}} ,\]
where
\[  \gamma:= \frac{\delta}{\sqrt{K}} + \sqrt{ \frac{\delta^2}{K} + n-1    }.\]
The minimum value is achieved when
\[  N= n +\frac{\delta\gamma}{\sqrt{K}}
  .\]
This gives (\ref{diamM}). In this case,
\[ H = \frac{K}{ \gamma^2}.\]

It is tricky to estimate the volume of the manifolds since there is no upper bound on ${\rm Vol}(B(p, R))$ only under the weighted Ricci curvature bound ${\rm Ric}^N \geq (N-1)H$.
We shall use (\ref{Volume**}) to get an upper bound on ${\rm Vol}(M) = {\rm Vol}(B_p(R))$ with $ R = \pi/\sqrt{H}$.

Let
\[ r:=  \frac{\pi}{2} \sqrt{\frac{n-1}{K}}.\]
Firstly suppose that
\[ d_p \leq r.\]
Then by (\ref{small3}), we get
\[ {\rm Vol} (M) ={\rm Vol}(B_p(r)) \leq \phi(p)K^{-n/2} C\Big ( n, \frac{\delta}{\sqrt{K}} \Big ),\]
where $C\Big ( n, \frac{\delta}{\sqrt{K}} \Big ):= (n-1)^{n/2} \int_{0}^{\pi/2} \sin^{n-1}(s) e^{\delta \sqrt{\frac{n-1}{K}} s} ds$.

We now assume that $ r < d_p \leq \pi\gamma/\sqrt{K}$. Let $\rho_o = r/2$. By Lemma \ref{lem4.3} (a), we have got two bounds  (\ref{small1}) and (\ref{small2}) on ${\rm Vol}(B_p(\rho_o))$ and ${\rm Vol} (B_p(r)\setminus B_p(\rho_o))$, respectively. On the other hand, from (\ref{RicH}) and by Lemma \ref{lem4.4}  and Theorem \ref{revol}, we can get an upper bound  on ${\rm Vol}(B_p(R))$ for $\chi(t) := [s_H(t)]^{N-1}$  and $R = \pi/\sqrt{H}$ by (\ref{Volume**}).
Then we obtain an upper bound on ${\rm Vol}(M)={\rm Vol}(B_p(R))$ as following
\begin{eqnarray*}
{\rm Vol}(M) & \leq & {\rm Vol}(B_p(\rho_o))+ \frac{ {\rm Vol} (B_p(r)\setminus B_p(\rho_o))}{\int_{\rho_o}^r \chi(t)dt }\int_{\rho_o}^R \chi(t) dt \\
& \leq & \phi(p) \omega_{n-1}  \Big ( \frac{n-1}{K}\Big )^{n/2} \Big \{ \int_0^{\pi/4} \sin^{n-1} (s) e^{\delta \sqrt{\frac{n-1}{K} } s } ds
 \\
&& +  \frac{\int_{\pi/4}^{\pi/2} \sin^{n-1}(s) e^{\delta \sqrt{\frac{n-1}{K} } s } ds}{ \int_{h}^{2h} \sin^{N-1}(s)ds }  \int_{h}^{\pi} \sin^{N-1} (s) ds\Big \},
\end{eqnarray*}
where
\[ h:= \frac{\pi}{4} \frac{ \sqrt{n-1}}{ \gamma}.\]
Thus there is a constant $C(n, \frac{\delta}{\sqrt{K}})$ depending only on $n$ and $\frac{\delta}{\sqrt{K}}$ such that
\[ {\rm Vol} (M) \leq \phi(p) K^{-n/2} C\Big ( n, \frac{\delta}{\sqrt{K} }\Big ).\]
Note that when  $\delta =0$, $ \gamma =\sqrt{n-1}$, $N=n $ and $H=K/(n-1)$. Then
\[{\rm Vol}(M) \leq \phi(p) \omega_{n-1}\Big ( \frac{n-1}{K}\Big )^{n/2} \int_0^{\pi}\sin^{n-1}(s) ds.\]
Clearly, this is a sharp bound. This completes the proof of Theorem \ref{prop4.5}. \qed

\section{Poincar\'{e}-Lichnerowicz inequality}\label{Sec3}

Let $(M, F, m)$ be a Finsler manifold equipped with a measure $m$ on $M$. In order to prove Theorem \ref{PLineq}, we need some necessary lemmas. Firstly, we have the following lemma.

\begin{lem}{\rm (\cite{Oh3})} \label{Lemma11.4} If $(M, F, m)$  satisfies that $\Lambda_{F}<\infty, \ {m}(M)< \infty$ and is complete, then the constant function $1$ belongs to $H_{0}^{1}(M)$, where $\Lambda_{F}$ denotes the reversibility constant of $F$.
\end{lem}

In our discussions, $(M, F, m)$ always satisfies that $\Lambda_{F}<\infty$.  Further,  we have the following result.

\begin{lem}\label{fsquareL}
For any $f, u \in H_{\mathrm{loc}}^{1}(M)$ satisfying $d f= 0$ almost everywhere on $M \backslash M_{u}$, we have
\be
\Delta^{\nabla u}f^{2}=2f \Delta^{\nabla u}f + 2 g_{\nabla u}\left(\nabla ^{\nabla u}f ~, \nabla ^{\nabla u}f\right), \label{eqLemma3.1}
\ee
equivalently,
\be
f \Delta^{\nabla u}f = \frac{1}{2}\Delta^{\nabla u}f^{2}- g_{\nabla u}\left(\nabla ^{\nabla u}f ~, \nabla ^{\nabla u}f\right). \label{squareLa}
\ee
\end{lem}

\noindent{\it Proof.} By (\ref{weigraLa}), we have
\[
\nabla^{\nabla u}f^{2} = 2f \nabla^{\nabla u}f.
\]
Further, by (\ref{DLOper}) and (\ref{weigraLa}), we have
\beqn
\Delta ^{\nabla u} f^{2}&=& {\rm div}_{m}\left(\nabla^{\nabla u}f^{2}\right)= {\rm div}_{m}(2 f \nabla^{\nabla u}f)\\
&=& 2f \Delta ^{\nabla u} f + 2 df (\nabla^{\nabla u}f).
\eeqn
It is easy to see that
\[
 df (\nabla^{\nabla u}f)=g_{\nabla u}\left(\nabla ^{\nabla u}f ~, \nabla ^{\nabla u}f\right).
\]
Thus we get (\ref{eqLemma3.1}). \qed

\vskip 2mm

In the following, we say that ${\rm Ric}^{\infty} \geq K$ for some $K \in {R}$ if ${\rm Ric}^{\infty}(v) \geq K F^{2}(x, v)$ for all $x\in M$ and $v\in T_{x}M$.
The following integrated form of the improved Bochner inequality is important for our discussion in this section.

\begin{lem}\label{IntformBo} {\rm (Integrated form,  \cite{Oh3})}  Assume that ${\rm Ric}^{\infty} \geq K$ for some $K \in {R}.$ Then, given $u \in H_{\mathrm{loc}}^{2}(M) \cap {\cal C}^{1}(M)$ such that $\Delta u \in H_{\mathrm{loc}}^{1}(M),$ we have
\beq
&&-\int_{M} d \phi\left(\nabla^{\nabla u}\left[\frac{F^{2}({\nabla} u)}{2}\right]\right) d{m} \nonumber\\
&& \geq \int_{M} \phi\left\{d({\Delta} u)({\nabla} u)+K F^{2}({\nabla}u)+d[F({\nabla}u)]\left(\nabla^{\nabla u}[F({\nabla} u)]\right)\right\} d{m} \label{Intform}
\eeq
for all nonnegative functions $\phi \in H_{0}^{1}(M) \cap L^{\infty}(M)$.
\end{lem}

The following lemma is derived from Lemma \ref{IntformBo} and is necessary for the proof of Theorem \ref{PLineq}.
\begin{lem}
Assume that $(M, F, {m})$ is compact and satisfies ${\rm Ric}^{\infty} \geq K>0$. Then , for any $u \in H^{2}(M) \cap \mathcal{C}^{1}(M)$ such that $\Delta u \in H^{1}(M)$, we have
\be
\int_{M} F^{2}({\nabla}u) d{m} \leq \frac{1}{K}\left\{ \int_{M}({\Delta} u)^{2} d{m} - \int_{M} g_{\nabla u}\left(\nabla ^{\nabla u}F(\nabla u) , \nabla ^{\nabla u}F(\nabla u)\right)dm \right\}. \label{mainlem}
\ee
\end{lem}

\noindent{\it Proof.} From (\ref{Intform}) and by taking test function $\phi \equiv 1$ according to Lemma \ref{Lemma11.4}, we have the following
\be
\int_{M} \left\{d({\Delta} u)({\nabla} u)+K F^{2}({\nabla}u)+d[F({\nabla}u)]\left(\nabla^{\nabla u}[F({\nabla} u)]\right)\right\} d{m} \leq 0.\label{firstequ}
\ee
Noticing that
\beqn
&& \int_{M} d({\Delta} u)({\nabla} u)dm = - \int_{M} (\Delta u)^{2}dm = -\|\Delta u \|^{2}_{L^{2}},\\
&& \int_{M} d[F({\nabla}u)]\left(\nabla^{\nabla u}[F({\nabla} u)]\right) = - \int_{M} F({\nabla}u) \Delta^{\nabla u}\left[F({\nabla}u)\right]dm ,
\eeqn
(\ref{firstequ}) becomes
\[
 K \int_{M}F^{2}({\nabla}u)dm \leq  \int_{M} (\Delta u)^{2}dm +   \int_{M} F({\nabla}u) \Delta^{\nabla u}\left[F({\nabla}u)\right]dm .
\]

Further, by Lemma \ref{fsquareL}, we have
\beqn
K \int_{M}F^{2}({\nabla}u)dm &\leq&  \int_{M} (\Delta u)^{2}dm +   \int_{M} \Delta ^{\nabla u}\left[\frac{ F^{2}({\nabla}u)}{2} \right]dm \\
    && -\int_{M} g_{\nabla u}\left(\nabla ^{\nabla u}F(\nabla u) , \nabla ^{\nabla u}F(\nabla u)\right)dm \\
   &=& \int_{M} (\Delta u)^{2}dm  -\int_{M} g_{\nabla u}\left(\nabla ^{\nabla u}F(\nabla u) , \nabla ^{\nabla u}F(\nabla u)\right)dm .
\eeqn
From this, we obtain (\ref{mainlem}). \qed

\vskip 2mm

In the following, normalizing $m$ as $m(M)=1,$ we define the variance of $f \in L^{2}(M)$ as
\be
\operatorname{Var}_{m}(f):=\int_{M}\left(f-\int_{M} f d m \right)^{2} d m = \int_{M} f^{2} d m -\left(\int_{M} f d m \right)^{2}.
\ee
Further, for any $f\in H^{1}(M)$, let $(u_{t})_{t\geq 0}$ be the global solution to the heat equation $\pa _{t}u_{t} = \Delta u_{t}$ with $u_{0}= f$. By the properties of the linearized heat semigroups, we have the mass conservation
\be
\int_{M} u_{t} dm = \int_{M} f dm
\ee
and ergodicity
\be
\lim\limits_{t\rightarrow \infty} u_{t} = \int_{M} f d{m} \ \ \ \label{eq15.2}
\ee
in $L^{2}(M)$ (see \cite{Oh2}\cite{Oh3}).

Now we are in the position to prove Theorem \ref{PLineq}.
\vskip 2mm

\noindent {\it Proof of Theorem \ref{PLineq}.} \ Let $(u_{t})_{t\geq 0}$ be the global solution to the heat equation with $u_{0}= f$. Put $\Phi (t):= \|u_{t}\|^{2}_{L^2}= \int_{M} u^{2}_{t} d m$. Then the ergodicity (\ref{eq15.2}) implies that
\[
\lim\limits_{t \rightarrow \infty}\Phi (t)= \lim\limits_{t \rightarrow \infty}\int_{M} u^{2}_{t} d m = \left(\int_{M} f d m \right)^{2}
\]
and
\beqn
{\rm Var}_{m}(f)&=&  \int_{M} f^{2} d m -\left(\int_{M} f d m \right)^{2}\\
                &=&  \Phi (0) - \lim\limits_{t \rightarrow \infty} \Phi (t)= - \int_{0}^{\infty} \Phi ' (t) dt .
\eeqn
On the other hand, by the definition of $\Phi (t)$ , we know that
\beq
\Phi^{\prime}(t)& =& 2 \int_{M} u_{t} \Delta u_{t} d{m} =  -2 \int_{M}du_{t}(\nabla u_{t}) d m \nonumber\\
                &=& -2 \int_{M} F^{2}\left(\boldsymbol{\nabla} u_{t}\right) d{m}.\label{Phiprime}
\eeq

By (4.2) in \cite{OS2}, we have the following
\be
\frac{\partial}{\partial t}\left[F^{2}\left({\nabla} u_{t}\right)\right]=2 d\left({\Delta} u_{t}\right)\left({\nabla} u_{t}\right) \label{Le14.1}
\ee
for all $t >0$. Hence, from (\ref{Phiprime}) and by (\ref{Le14.1}), we have
\[
\Phi^{\prime \prime}(t)=-4 \int_{M} d\left(\Delta u_{t}\right)\left(\nabla u_{t}\right) d {m}=4 \int_{M}\left(\Delta u_{t}\right)^{2} dm =4\left\|\Delta u_{t}\right\|_{L^{2}}^{2}
\]
for all $t>0$. Thus, from (\ref{mainlem}) for $u_{t}$ and by (\ref{Phiprime}), we get the following
\[
-\frac{1}{2}\Phi^{\prime}(t) \leq \frac{1}{4K} \Phi^{\prime \prime}(t)-\frac{1}{K}\int_{M}g_{\nabla u_{t}}\left(\nabla ^{{\nabla u_{t}}}F({\nabla u_{t}}), \nabla ^{{\nabla u_{t}}}F({\nabla u_{t}})\right)dm,
\]
that is,
\[
- \Phi^{\prime}(t) \leq \frac{1}{2K} \Phi^{\prime \prime}(t)-\frac{2}{K}\int_{M}g_{\nabla u_{t}}\left(\nabla ^{{\nabla u_{t}}}F({\nabla u_{t}}), \nabla ^{{\nabla u_{t}}}F({\nabla u_{t}})\right)dm.
\]
Then,
\beqn
&& {\rm Var}_{m}(f)=-\int_{0}^{\infty} \Phi ^{\prime}(t)dt \\
&& \leq \frac{1}{2K}\left(\lim\limits_{t \rightarrow \infty}\Phi ^{\prime}(t) - \Phi ^{\prime}(0)\right)-\frac{2}{K}\int_{M}\left(\int_{0}^{\infty} g(t)dt\right)dm,
\eeqn
where
\[
g(t):= g_{\nabla u_{t}}\left(\nabla ^{{\nabla u_{t}}}F({\nabla u_{t}}), \nabla ^{{\nabla u_{t}}}F({\nabla u_{t}})\right).
\]
From (\ref{Phiprime}), $\Phi ^{\prime}(0)= -2 \int_{M}F^{2}(\nabla f) d m$. Further, it is not difficult to prove that $\lim\limits_{t\rightarrow \infty}{\cal E}(u_{t})=0$ (e.g. see \cite{Oh2}\cite{OS1}), and then, $\lim\limits_{t\rightarrow \infty} \Phi^{\prime}(t) =0$.  Thus, we get the following
\be
{\rm Var}_{m}(f)\leq \frac{1}{K}\int_{M} F^{2}(\nabla f) dm - \frac{2}{K}\int_{M}\left(\int_{0}^{\infty} g(t)dt\right)dm. \label{PLine2}
\ee
This completes the proof of Theorem \ref{PLineq}. \qed

\vskip 2mm

As an application of Theorem \ref{PLineq}, the following result is natural.

\begin{cor} Suppose that $(M, F, m)$ is closed and satisfies $m(M)=1$ and $\operatorname{Ric}^{\infty} \geq K>0$.  Let $\lambda_1$ be the first (nonzero) eigenvalue of the Finsler Laplacian. Then
\be
\lambda_1 \geq K+\delta ,
\ee
where
\[ \delta := \inf_{f} \frac{2}{{\rm Var}_m(f)}\int_{M}\left(\int_{0}^{\infty} g(t)dt\right)dm.\]

\end{cor}

\vskip 8mm

\noindent
Xinyue Cheng \\
School of Mathematical Sciences \\
Chongqing Normal University \\
Chongqing  401331,  P. R. of China  \\
E-mail: chengxy@cqnu.edu.cn

\vskip 3mm

\noindent
Zhongmin Shen \\
Department of Mathematical Sciences \\
Indiana University-Purdue University Indianapolis \\
IN 46202-3216, USA  \\
E-mail: zshen@math.iupui.edu

\end{document}